\documentclass[12pt]{article}
\usepackage{amsmath,amssymb,amscd,verbatim,amsthm}
\usepackage{latexsym, epsfig, graphpap, graphics, mathrsfs}
\usepackage[varg]{pxfonts}
\usepackage{mathrsfs}

\begin{document}
 \noindent{\bf{Some new lacunary $f$-statistical $A$-convergent sequence spaces of order $\alpha$ }}

\vskip 0.5 cm

\noindent Ekrem Sava$\c{s}$

\noindent Istanbul Commerce University, 34840 Istanbul, Turkey

\noindent E-mail : ekremsavas@yahoo.com

\vskip 0.5 cm

\noindent Stuti Borgohain $^*${\footnote {The work of the
authors was carried under the Post Doctoral Fellow under National Board of Higher 
Mathematics, DAE, project No. NBHM/PDF.50/2011/64}}

\noindent Department of Mathematics

\noindent Indian Institute of Technology, Bombay

\noindent Powai:400076, Mumbai, Maharashtra; INDIA.
\vskip 0.5 cm

\noindent E-mail : stutiborgohain$@$yahoo.com
\vskip 1 cm

\noindent{\footnotesize {\bf{Abstract:}} We study the concept of density for sets of natural numbers in some lacunary $A$-convergent sequence spaces. Also we are trying to investigate some relation between the ordinary convergence and module statistical convergence for evey unbounded modulus function. Morever we also study some results on the newly defined lacunary $f$-statistically $A$-convergent sequence spaces with respect to some Musielak-Orlicz function.} \\

\section {Introduction}

In order to extend the notion of convergence of sequences,  Fast \cite{Fast} and  Schoenberg \cite{Schoenberg} independently introduced the concept of statistical convergence. Later on it was studied from sequence space point of view and also linked with summability theory by many mathematicians including  $\check{S}$al$\grave{a}$t \cite{Salat}, Fridy \cite{Fridy}, Tripathy and Borgohain \cite{Tripathy} and many more. Kolk \cite{Kolk} began to study the applications of statistical convergence to Banach spaces. Connor et. al \cite{Connor} proved some important results that relate the statistical convergence to classical properties of Banach spacs.\\

The notion depends on the idea of asymtotic density of subsets of the set $\mathbb{N}$ of natural numbers. 
A subset $A$ of $\mathbb{N}$ is said to have natural density ${\delta(A)}$ if

$${\delta(A)}={\lim_{n \rightarrow \infty}}{\frac{1}{n}}{\sum\limits_{k=1}^n}{\chi_A}(k).$$

\noindent where $\chi_A$ is the characteristic function of $A$.\\

A sequence  $(x_n)_n$ is said to be statistically convergent to $L$, if for any $\varepsilon > 0$,we have $\delta\{(k \in \mathbb{N}:\vert x_k - L \vert \geq \varepsilon )\} = 0.$ Analogously, $(x_n)_n$ is said to be statistically Cauchy if for each $\varepsilon>0$ and $n \in \mathbb{N}$ there exists an integer $m \geq n$ such that $d(\{i \in \mathbb{N}:\Vert x_i - x_m \Vert < \varepsilon \})=1.$\\

Fridy \cite{Fridy} proved that in a Banach space, a sequence is statistically convergent if and only if it is statistically Cauchy. Fast \cite{Fast} proved that st$\displaystyle\lim_n x_n =x$ if and only if there exists $A \subset \mathbb{N}$ with $d(A)=1$ and $\displaystyle\lim_{n \in A} x_n =x$.\\

The notion of a modulus function was introduced by Nakano \cite{Nakano}. Ruckle \cite{Ruckle} and Maddox \cite{Maddox} have introduced and discussed some properties of sequence spaces defined by using a modulus function. By the definition of modulus function, we mean a function $f: R^+ \rightarrow R^+ $ which satisfies:

\begin{enumerate}

\item $f(x)=0$ if and only if $x=0$.

\item $f(x+y) \leq f(x)+f(y)$ for every $x,y \in R^+$.

\item $f$ is increasing.

\item $f$ is continuous from the right at 0.

\end{enumerate}

A lacunary sequence is defined as an increasing integer sequence $\theta = (k_r)$ such that $k_0=0$ and $h_r=k_r-k_{r-1} \rightarrow \infty$ as $r \rightarrow \infty.$  Throughout this paper, the intervals determined by $\theta$ will be denoted by $J_r=(k_{r-1}, k_r]$ and the ratio $\frac{k_r}{k_{r-1}}$ will be defined by $\phi_r$. \\

In this paper, we study on a new lacunary $f$-statistical $A$-convergent sequence space of order $\alpha$ with respect to the Musielak-Orlicz fuction. We also investigate some results on a new concept of nonmatrix convergence which is intermediate between the ordinary convergence and the statistical convergence.

\section{Definitions and basic results}

By the definition of $f$-density of a set $A \subseteq \mathbb{N}$ , we mean $d_f(A)=\displaystyle\lim_n \frac{f(\vert A(n) \vert)}{f(n)}$, (in case this limit exists )where $f$ is an unbounded modulus function. \\

Let $X$ be a normed space and let $(x_n)_n$ be a sequence in $X$. If if for each $\varepsilon > 0$, $d_f(\{i \in \mathbb{N}: \Vert x_i-x \Vert > \varepsilon \})=0$, then it is said that the $f$-statistical limit of $(x_n)_n$ is $x \in X$, and we write it as $f$-stlim$x_n=x$. \\

It is clear that $d(A)=1-d(\mathbb{N} \backslash A)$.\\

Let us assume that $A \subseteq \mathbb{N}$ and $d_f(A)=0$. For every $n \in N$ we have $f(n) \leq f(\vert A(n)\vert)+f(\vert (\mathbb{N} \backslash A)(n)\vert )$ and so

$$1 \leq \frac{f(\vert A(n) \vert)}{f(n)} +\frac{f(\vert ( \mathbb{N} \backslash A)(n)\vert )}{f(n)} \leq \frac{f(\vert A(n) \vert)}{f(n)}+1$$

By taking limits we deduce that $d_f(\mathbb{N} \backslash A)=1$.\\

{\bf{Corollary 2.1. }}Let $f, g$ be unbounded moduli, $X$ a normed space, $(x_n)_n$ a sequence in $X$ and $x, y \in X$. We have

\begin{enumerate}

\item The $f$-statistical convergence implies the statistical convergence (to the same limit).

\item The $f$-statistical limit is unique whenever it exists.

\item Moreover, two different methods of statistical convergence are always compatible, which means that if $f$-st$\displaystyle\lim x_n =x$ and $g$-st$\displaystyle\lim x_n =y$ then $x=y$.

\end{enumerate}

By an Orlicz function , we mean a function $M: [0,\infty )\rightarrow [0,\infty )$, which is continuous, non-decreasing and convex with $M(0) = 0, M(x)>0$, for $x>0$ and $M(x)\rightarrow \infty$, as $x\rightarrow \infty$.\\

The idea of Orlicz function is used to construct the sequence space, (see Lindenstrauss and Tzafriri \cite{Lindenstrauss}),

$$\ell_M=\left\{ (x_k) \in w: \displaystyle\sum_{k=1}^\infty M \left(\frac{\vert x \vert}{\rho}\right) < \infty, \mbox{~for some~} \rho>0 \right\}. $$

This space $\ell_M$ with the norm,

$$ \Vert x \Vert = \mbox{inf}\left\{ \rho>0: \displaystyle\sum_{k=1}^\infty M \left(\frac{\vert x_k \vert}{\rho}\right) \leq 1 \right\}$$

becomes a Banach space which is called an Orlicz sequnce space.\\

Musielak \cite{Musielak} defined the concept of Musielak-Orlicz function as $\mathscr{M}=(M_k)$.\\
A sequence $\mathscr{N}=(N_k)$ defined by 
$$N_k(v)=\sup \{ \vert v \vert u -M_k(u): u \geq 0 \}, k=1,2,..$$

is called the complementary function of a Musielak-Orlicz function $\mathscr{M}$. The Musielak-Orlicz sequence space $t_\mathscr{M}$ and its subspace $h_\mathscr{M}$ are defined as follows:

$$t_\mathscr{M}=\{ x \in w: I_\mathscr{M}(cx) < \infty \mbox{~for some ~} c>0\},$$
$$h_\mathscr{M}=\{x \in w:I_\mathscr{M}(cx) < \infty, \forall c >0\},$$

where $I_\mathscr{M}$ is a convex modular defined by,

$$I_\mathscr{M}(x) = \displaystyle\sum_{k=1}^\infty M_k(x_k), x=(x_k) \in t_\mathscr{M}.$$

It is considered $t_\mathscr{M}$ equipped with the Luxemberg norm

$$\Vert x \Vert = \inf \left\{ k>0: I_\mathscr{M}\left(\frac{x}{k}\right) \leq 1 \right\}$$

or equiped with the Orlicz norm

$$\Vert x \Vert^0 = \inf \left\{ \frac{1}{k} (1+I_\mathscr{M}(kx)): k>0 \right\}.$$

A Musielak-Orlicz function $(M_k)$ is said to satisfy $\Delta_2$-condition if there exist constants $a, K>0$ and a sequence $c=(c_k)_{k=1}^\infty \in \ell_+^1$ ( the positive cone of $\ell^1$) such that the inequality

$$M_k(2u) \leq KM_k(u)+c_k$$

holds for all $k \in N$ and $u \in R_+$, whenever $M_k(u) \leq a$.

For any lacunary sequence $\theta= (k_r)$, the space $N_\theta$ defined as, (Freedman et al.[1])
$$N_\theta=\left\{(x_k): \displaystyle\lim_{r \rightarrow \infty} h_r^{-1} \displaystyle\sum_{k \in J_r} \vert x_k -L \vert =0, \mbox{~for some~} L \right\}.$$

The space $N_\theta$ is a $BK$ space with the norm,
$$\Vert (x_k) \Vert_\theta= \displaystyle\sup_r h_r^{-1} \displaystyle\sum_{k \in J_r} \vert x_k \vert.$$

A sequence $(x_i)$ is said to be lacunary $f$-statistically $A$-convergent to $L$ if,
$$d_f\left(\left\{ k \in N: \displaystyle\lim_{r \rightarrow \infty} \frac{1}{h_r^\alpha} \displaystyle\sum_{i \in I_r}  M_i \left(\frac{\vert A_i(x) -L \vert}{\rho^{(k)}} \right)  > \varepsilon, \mbox{~for some~} L \mbox{~and~} \rho^{(k)}>0 \right\} \right) =0,$$

where $I_r = (i_{i-1}, i_r]$ , $A=(A_i(x))$ such that $A_i x =\displaystyle\sum_{k=1}^\infty a_{ik} x_k$ converges for each $i$ , $h_r^\alpha$ denotes the $\alpha$-th power $(h_r^\alpha)$ of $h_r$, that is, $h^\alpha=(h_r^\alpha)=(h_1^\alpha, h_2^\alpha,...h_r^\alpha,...)$ and $\mathscr{M}=(M_k)$ is a Musielak-Orlicz function and it  is written as $f_\theta^\alpha(A, \mathscr{M})$-statistically convergentto $L$ with respect to the Musielak-Orlicz function $\mathscr{M}$.\\

Some particular cases : \\

If we take $\alpha=1$, then the lacunary $f$-statistically $A$-convergence of order $\alpha$ reduces to the lacunary $f$-statistically $A$-convergent to $L$ ($f_\theta(A, \mathscr{M})$-statistically convergent to $L$) , i.e.
$$d_f\left(\left\{ k \in N: \displaystyle\lim_{r \rightarrow \infty} \frac{1}{h_r} \displaystyle\sum_{i \in I_r}  M_i \left(\frac{\vert A_i(x) -L \vert}{\rho^{(k)}} \right)  > \varepsilon, \mbox{~for some~} L \mbox{~and~} \rho^{(k)}>0 \right\} \right) =0.$$\\
 
If $\theta=(2^r)$ and $\alpha=1$ , then the sequence $(x_i)$ is said to be $f(A,\mathscr{M})$-statistically convergent to $L$ if ,
$$d_f\left(\left\{ k \in N:  \displaystyle\sum_{i \in I_r}  M_i \left(\frac{\vert A_i(x) -L \vert}{\rho^{(k)}} \right)  > \varepsilon, \mbox{~for some~} L \mbox{~and~} \rho^{(k)}>0 \right\} \right) =0.$$

If $M_i(x)=x$, $\theta=(2^r)$ and $\alpha=1$,  then we can say the sequence $(x_i)$ is $f$-statistically $A$-convergent to $L$ , if
$$d_f\left(\left\{ k \in N:  \displaystyle\sum_{i \in I_r}  \vert A_i(x) -L \vert   > \varepsilon, \mbox{~for some~} L \right\}\right)=0.$$

Let us define a space of lacunary $A$-convergent sequences of order $\alpha$ defined by Musielak-Orlicz function as,
$$w_\theta^\alpha(A, \mathscr{M})=\left\{ (x_i):\displaystyle\lim_{r \rightarrow \infty} \frac{1}{h_r^\alpha} \displaystyle\sum_{i \in I_r} M_i \left(\frac{\vert A_i(x) - L \vert}{\rho^{(i)}} \right)=0, \mbox{~for some~} L \mbox{~and~} \rho^{(i)}>0\right\}.$$

In particular, if we have  $M_i(x)=x$, $\theta=(2^r)$ and $\alpha=1$, then $w_\theta^\alpha(A, \mathscr{M})$ reduces to,
$$w(A)=\left\{ (x_i): \displaystyle\sum_{i \in I_r} \vert A_i(x) -L \vert =0, \mbox{~for some~} L \right\}.$$

\section{Main Results}

{\bf{Theorem 3.1.}}  Let $(x_n)_n$ be a sequence in a normed space $w_\theta^\alpha(A, \mathscr{M})$ and $f$ an unbounded modulus. Then $f$-st$\displaystyle\lim x_i=L$ if and only if there exists $X \subseteq \mathbb{N}$ such that $d_f(X)=0$ and $\displaystyle\lim_{i \in \mathbb{N} \backslash X} x_i =L $.\\

{\it{Proof}}: Let $B_j = \left\{ i \in \mathbb{N}:\displaystyle\lim_{r \rightarrow \infty} \frac{1}{h_r^\alpha} \displaystyle\sum_{i \in I_r} M_i \left(\frac{ \vert A_i(x)- L \vert}{\rho^{(i)}} \right)  > \frac{1}{j} \right\}$, for every $j \in \mathbb{N}$.\\

Since $B_j \subset B_{j+1}$ and $d_f(B_j)=0$, it is required to prove the case where some of the $B_j$'s are non-empty, in particular we can assume that $B_1 \neq \phi$. \\

Choose any $r_1 \in B_1$. Now, by taking $r_2 \in B_2$ with $r_2 > r_1$ and $\frac{f(\vert B_2(i) \vert)}{f(i)} \leq \frac{1}{2}$, if $i \geq r_2$.\\

Inductively we obtain $r_1 < r_2 < r_3 ...$ such that $r_j \in B_j$ and $\frac{f(\vert B_j(i) \vert)}{f(i)} \leq \frac{1}{j}$ whenever $i \geq r_j$.\\

Now consider $X= \displaystyle\cup_{j \in \mathbb{N}} ([r_j, r_{j+1}) \cap B_j)$. Then for every $i \geq r_1$ there exists $j \in \mathbb{N}$ such that $r_j \leq i \leq r_{j+1}$ and if $n \in X(i)$ then $n < r_{j+1}$, which implies $n \in B_j$. Therefore $X(i) \subseteq B_j(i)$ and thus,

$$\frac{f(\vert X(i) \vert )}{f(i)} \leq \frac{f(\vert B_j(i) \vert )}{f(i)} \leq \frac{1}{j}$$

which follows $d_f(X)=0$.\\

For $\varepsilon > 0$ and $j \in \mathbb{N}$ such that $\frac{1}{j} < \varepsilon$, we have $i \in \mathbb{N} \backslash X$ and $i \geq r_j$ for which  there exists $k \geq j$ with $r_k \leq i \leq r_{k+1}$ and this implies $i \notin B_k$, so,

$$\displaystyle\lim_{r \rightarrow \infty} \frac{1}{h_r^\alpha} \displaystyle\sum_{i \in I_r} M_i \left( \frac{\vert A_i(x) - L \vert}{\rho^{(i)}} \right) \leq \frac{1}{k} \leq \frac{1}{j} < \varepsilon.$$ 

We conclude $\displaystyle\lim_{i \in \mathbb{N} \backslash X} x_i =L$.\\

Conversely, let us assume that $X \in \mathbb{N}$ satisfies $\displaystyle\lim_{i \in \mathbb{N} \backslash X} x_i=L$ and $d_f(X)=0$. For $\varepsilon > 0$, there exists $i_0 \in \mathbb{N}$ such that if $i > i_0$ and $i \in \mathbb{N} \backslash X$ then 

$$\displaystyle\lim_{r \rightarrow \infty} \frac{1}{h_r^\alpha} \displaystyle\sum_{i \in I_r} M_i \left( \frac{\vert A_i(x) - L \vert}{\rho^{(i)}} \right) \leq \varepsilon.$$

This implies,

$$\left\{i \in \mathbb{N}: \displaystyle\lim_{r \rightarrow \infty} \frac{1}{h_r^\alpha} \displaystyle\sum_{i \in I_r} M_i \left( \frac{\vert A_i(x) - L \vert}{\rho^{(i)}} \right) > \varepsilon \right\} \subseteq X \cup \{1,..i\}$$

and then $d_f\left(\left\{i \in \mathbb{N}: \displaystyle\lim_{r \rightarrow \infty} \frac{1}{h_r^\alpha} \displaystyle\sum_{i \in I_r} M_i \left( \frac{\vert A_i(x) - L \vert}{\rho^{(i)}} \right) > \varepsilon \right\}\right)=0$.\\

\vskip 0.3 cm

{\bf{Definition 3.2.}} The sequence $(x_n)_n$ is $f$-statistically Cauchy if for every $\varepsilon > 0$ there exists $N \in \mathbb{N}$ such that $d_f(\{i \in \mathbb{N} : \Vert x_i - x_N \Vert > \varepsilon \})= 0$.\\

It is well known that if a sequence is $f$-statistically convergent then it is $f$-statistically Cauchy. The converge is true if the space is complete and this result is a particular case of filter convergence.\\

{\bf{Thoerem 3.3.}}  For the Banach space $w_\theta^\alpha(A,\mathscr{M})$ where $f$ be an unbounded modulus and $(x_n)_n$ is an $f$-statistically Cauchy sequence. Then $(x_n)_n$ is $f$-statistically convergent.\\

{\it{Proof}}: For every $k \in N$, let $x_{N_k}$ be such that 

$$d_f\left(\left\{i \in \mathbb{N}: \displaystyle\lim_{r \rightarrow \infty} \frac{1}{h_r^\alpha} \displaystyle\sum_{i \in I_r} M_i \left( \frac{\vert A_i(x) - A_{N_k}(x) \vert}{\rho^{(i)}} \right)> \frac{1}{k} \right\}\right)=0.$$

Consider the sets $I_k = \cap_{j \leq k} B(x_{N_j}, \frac{1}{j})$ and $J_k=\{ i \in \mathbb{N}: x_i \notin I_k\}.$\\

Then for each $k \in \mathbb{N}$ we have $\mbox{diam}(I_k) \leq \frac{2}{k}$ and 

$$J_k=\cup_{j \leq k} \left\{i \in \mathbb{N}: \displaystyle\lim_{r \rightarrow \infty} \frac{1}{h_r^\alpha} \displaystyle\sum_{i \in I_r} M_i \left( \frac{\vert A_i(x) - A_{N_j}(x) \vert}{\rho^{(i)}} \right) > \frac{1}{j} \right\},$$

which implies $d_f(J_k)=0$.\\

Consequently each $I_k \neq \phi$.\\

Sincet $I_1 \supseteq I_2 \supseteq ...$ and $J_1 \subseteq J_2 \subseteq...$ so as in the proof of the previous theorem, we can find a sequence of natural numbers $r_1<r_2<....$ such that $r_k \in J_k$ and if $i \geq r_k$ then $\frac{f(\vert J_k(i)\vert )}{f(i)} \leq \frac{1}{k}.$ \\

Considering $X= \cup_{k \in \mathbb{N}} ([r_k, r_{k+1}) \cap J_k)$, we get $d_f(X)=0$.\\

Since $w_\theta^\alpha(A,\mathscr{M})$ is complete, then $\cap_{k \in \mathbb{N}} I_k$ has exactly one element, say $x$.\\

We have to prove that $\displaystyle\lim_{i \in \mathbb{N} \backslash X} x_i =L$.\\

For $\varepsilon > 0$, choose $j \in \mathbb{N}$ with $\frac{2}{j} < \varepsilon$. If $i \geq r_j$ and $i \in \mathbb{N} \backslash X$ then there exists $k \geq j$ such that $r_k \leq i < r_{k+1}$ and then $i \notin J_k$, which implies $x_i \in I_k$ and thus, 

$$ \displaystyle\lim_{r \rightarrow \infty} \frac{1}{h_r^\alpha} \displaystyle\sum_{i \in I_r} M_i \left( \frac{\vert A_i(x) - A_{N_k}(x) \vert}{\rho^{(i)}} \right) \leq \frac{2}{k} \leq \frac{2}{j} < \varepsilon.$$

This completes the proof of the theorem.\\

{\bf{Theorem 3.4.}}  Let $(x_n)_n$ be a sequence in $w_\theta^\alpha(A,\mathscr{M})$. If for every unbounded modulus $f$ there exists $f$-st$\displaystyle\lim x_i$ then all these limits are the same $x \in w_\theta^\alpha (A, \mathscr{M})$ and $(x_n)_n$ also converges to $L$ in the norm topology.\\

{\it{Proof}} : It is proved that for $f,g$ two unbounded moduli, $X$ a normed space, $(x_n)_n$ a sequence in $X$ and $L, M \in X$, the $f$-statistical limit is unique whenever it exists. \\

If it is false that $\displaystyle\lim x_i =L$, there exists $\varepsilon > 0$ such that

$$X=\left\{ i \in \mathbb{N}: \displaystyle\lim_{r \rightarrow \infty} \frac{1}{h_r^\alpha} \displaystyle\sum_{i \in I_r} M_i \left( \frac{\vert A_i(x) - L \vert}{\rho^{(i)}} \right)> \varepsilon \right\}$$

is infinite.\\

Now, by choosing an unbounded modulus $f$ which will satisfy $d_f(X)=1$, then this clearly contradicts the assumption that $f$-stlim $x_i=L$, which completes the proof of the Theorem.\\

\textbf{Theorem 3.5.}  Let $\mathscr{M}$ be a Musielak-Orlicz function, $x=(x_i)$ be a bounded sequence and $\theta=(i_r)$ be a lacunary sequence. If $\displaystyle\lim_{r \rightarrow \infty} \frac{h_r}{h_r^\alpha}=1$,then $x \in f_\theta^\alpha(A, \mathscr{M}) \Rightarrow x \in w_\theta^\alpha (A, \mathscr{M})$.\\

{\it{Proof}}: Suppose that $x=(x_i)$ be a bounded sequence and $ f_\theta^\alpha(A,\mathscr{M})-\mbox{lim} x_i =L$.\\

Then we have, 
$$d_f \left( \left\{ i \in \mathbb{N}: \frac{1}{h_r^\alpha} \displaystyle\sum_{i \in I_r} M_i \left(\frac{\vert A_i(x) - L \vert}{\rho^{(i)}}\right) > \varepsilon \right\}\right) =0. $$

For $\varepsilon >0$  given, let us denote $\sum_1$ as the sum over $i \in I_r$, $\vert A_i(x)-L \vert \geq \varepsilon$ and $\sum_2$ denotes the sum over $i \in I_r$, $\vert A_i(x)-L \vert < \varepsilon$ respectively.\\ 

As $x \in \ell_\infty$, we have a constant $T>0$ such that $\vert x_i \vert \leq T$. Now, for $\varepsilon > 0$, we have,\\

$\frac{1}{h_r^\alpha} \displaystyle\sum_{i \in I_r} M_i \left(\frac{\vert A_i(x) - L \vert}{\rho^{(i)}}\right)$\\

$\frac{1}{h_r^\alpha} \sum_1 M_i\left(\frac{\vert A_i(x)-L \vert}{\rho^{(i)}}\right) +
\frac{1}{h_r^\alpha} \sum_2 M_i\left(\frac{\vert A_i(x)-L \vert}{\rho^{(i)}}\right)$\\

$\leq \frac{1}{h_r^\alpha} \sum_1 \mbox{max} \left\{M_i\left(\frac{T}{\rho^{(i)}}\right), M_i\left(\frac{T}{\rho^{(i)}}\right) \right\} + \frac{1}{h_r^\alpha} \sum_2 M_i\left(\frac{\varepsilon}{\rho^{(i)}} \right)$\\

$\leq \mbox{max}\{M_i(K),M_i(K) \} \frac{1}{h_r^\alpha} \vert \{ i \in I_r: \vert A_i(x)-L \vert \geq \varepsilon \} \vert + \frac{h_r}{h_r^\alpha}  M_i(\varepsilon_1) , ~~\frac{T}{\rho^{(i)}}=K, \frac{\varepsilon}{\rho^{(i)}}=\varepsilon_1.$\\

Hence, $x \in w_\theta^\alpha(A, \mathscr{M})$. \\

{\bf{Theorem 3.6. }} Let $\mathscr{M}=(M_i)$ be a Musielak-Orlicz function where $(M_i)$ is pointwise convergent. Then, $w_\theta^\alpha(A, \mathscr{M}) \subset f_\theta^\alpha(A, \mathscr{M})$ if and only if $\displaystyle\lim_i M_i \left(\frac{ \nu}{\rho^{(i)}}\right) > 0$ for some $\nu>0, \rho^{(i)}>0$.\\

{\it{Proof}} : Let $\varepsilon >0$ and $x \in w_\theta^\alpha(A, \mathscr{M})$. \\

If $\displaystyle\lim_i M_i\left( \frac{\nu}{\rho^{(i)}}\right) >0$, then we have a number $c > 0$ such that

$$M_i \left(\frac{\nu}{\rho^{(i)}}\right) \geq c, \mbox{~for~} \nu > \varepsilon.$$

Let us define, $I_r^1=\left\{ i \in I_r:  M_i \left(\frac{\vert A_i(x)-L \vert}{\rho^{(i)}}\right) \geq \varepsilon \right\}$. \\

Then, 
\begin{eqnarray*}
\frac{1}{h_r^\alpha} \displaystyle\sum_{i \in I_r}  M_i \left(\frac{\vert A_i(x)-L  \vert}{\rho^{(i)}} \right)
&\geq &
 \frac{1}{h_r^\alpha} \displaystyle\sum_{i \in I_r^1}  M_i \left(\frac{\vert A_i(x)-L \vert}{\rho^{(i)}} \right)\\
&\geq &
c \frac{1}{h_r^\alpha} \vert A_0(\varepsilon) \vert
\end{eqnarray*}

Hence, it follows that $ x \in f_\theta^\alpha(A, \mathscr{M})$.\\

Conversely, let the condition do not exist. Now, for a number $\nu>0$ , let $\displaystyle\lim_i M_i\left(\frac{\nu}{\rho^{(i)}}\right)=0$ for some $\rho>0$. 

Select a lacunary sequence $\theta=(n_r)$ such that $M_i \left(\frac{\nu}{\rho^{(i)}}\right)< 2^{-r}$ for any $i > n_r$.\\

Also, if $A=I$ , then we can define a sequence $x$ by putting,

\[ A_i(x) = \left\{ \begin{array}{ll}
\nu & \mbox{if $n_{r-1} < i \leq \frac{n_r+n_{r-1}}{2}$};\\
\theta & \mbox{if $\frac{n_r+n_{r-1}}{2} < i \leq n_r$}.\end{array} \right. \]

Therefore,

\begin{eqnarray*}
\frac{1}{h_r^\alpha}\displaystyle\sum_{i \in I_r} M_i \left( \frac{\vert A_i(x) \vert}{\rho^{(i)}} \right)
&=&
\frac{1}{h_r^\alpha} \displaystyle\sum_{n_{r-1} < i \leq \frac{(n_r+n_{r-1})}{2}} M_i \left( \frac{\nu}{\rho^{(i)}} \right)\\
&<&
\frac{1}{h_r^\alpha} \frac{1}{2^{r-1}} \left[ \frac{n_r+n_{r-1}}{2} - n_{r-1} \right]\\
&=&
\frac{1}{2^r} \rightarrow 0 \mbox{~as~} r \rightarrow \infty.
\end{eqnarray*}

Thus we have $x \in w_\theta^{\alpha 0} (A, \mathscr{M})$. \\

But, 
\begin{eqnarray*}
\displaystyle\lim_{r \rightarrow \infty} \frac{1}{h_r^\alpha} \left\vert \left\{ i \in  I_r: \displaystyle\sum_{i \in I_r} M_i \left(\frac{\vert A_i(x) \vert}{\rho^{(i)}} \right) \geq \varepsilon \right\} \right\vert
&=&
\displaystyle\lim_{r \rightarrow \infty} \frac{1}{h_r^\alpha} \left \vert \left \{ i \in \left( n_{r-1}, \frac{n_r + n_{r-1}}{2} \right) : \displaystyle\sum_{i \in I_r} M_i \left( \frac{\nu}{\rho^{(i)}} \right) \geq \varepsilon \right\} \right\vert\\
&=&
\displaystyle\lim_{r \rightarrow \infty} \frac{1}{h_r^\alpha} \frac{n_r-n_{r-1}}{2} \\
&=& 
\frac{1}{2}.
\end{eqnarray*}

So, $x \notin f_\theta^\alpha(A, \mathscr{M})$.  \\

{\bf{Theorem 3.7. }} Let $\mathscr{M}=(M_i)$ be a Musielak-Orlicz function. Then $f_\theta^\alpha(A, \mathscr{M}) \subset w_\theta^\alpha(A, \mathscr{M})$ if and only if $\displaystyle\sup_\nu \displaystyle\sup_i M_i \left(\frac{\nu}{\rho^{(i)}}\right) < \infty$.\\

Proof: Let $x \in f_\theta^\alpha(A, \mathscr{M})$. Suppose $h(\nu)=\displaystyle\sup_i M_i \left(\frac{\nu}{\rho^{(i)}} \right)$ and $h=\displaystyle\sup_\nu h(\nu)$. Let

$$I_r^2=\left\{i \in I_r: M_i \left(\frac{\vert A_i(x)-L \vert}{\rho^{(i)}} \right) < \varepsilon \right\}.$$

Now, $M_i(\nu) \leq h$ for all $i, \nu >0$. So,
\begin{eqnarray*}
\frac{1}{h_r^\alpha} \displaystyle\sum_{i \in I_r}  M_i \left(\frac{\vert A_i(x)-L \vert}{\rho^{(i)}} \right)
&= &
 \frac{1}{h_r^\alpha} \displaystyle\sum_{i \in I_r^1}  M_i \left(\frac{\vert A_i(x)-L \vert}{\rho^{(i)}} \right)\\
&+&
 \frac{1}{h_r^\alpha} \displaystyle\sum_{i \in I_r^2} M_i \left(\frac{\vert A_i(x)-L \vert}{\rho^{(i)}} \right)\\
&\leq &
h \frac{1}{h_r^\alpha} \vert A_0(\varepsilon) \vert + h(\varepsilon).
\end{eqnarray*}

Hence, as $\varepsilon \rightarrow 0$, it follows that $x \in w_\theta^\alpha(A, \mathscr{M})$.\\

Conversely, suppose that
$$\displaystyle\sup_\nu \displaystyle\sup_i M_i \left( \frac{\nu}{\rho^{(i)}} \right) =\infty.$$

Then, we have
$$0< \nu_1< \nu_2 <...< \nu_{r-1}< \nu_r<...$$

so that $M_{n_r} \left(\frac{\nu_r}{\rho^{(i)}} \right) \geq h_r^\alpha$ for $r \geq 1$. For $A=I$, we set a sequence $x=(x_i)$ by,
\[ A_i(x) = \left\{ \begin{array}{ll}
\nu_r & \mbox{if $i=n_r$ for some $r=1,2,..$};\\
\theta & \mbox{otherwise}.\end{array} \right. \]

Then,

\begin{eqnarray*}
\displaystyle\lim_{r \rightarrow \infty} \frac{1}{h_r^\alpha} \left\vert \left\{ i \in I_r :  \displaystyle\sum_{i \in I_r} M_i \left(\frac{\vert A_i(x) \vert}{\rho^{(i)}} \right)  \geq \varepsilon \right\} \right \vert 
&=&
\displaystyle\lim_{r \rightarrow \infty} \frac{1}{h_r^\alpha} \\
&=&
0
\end{eqnarray*}

Hence, $x \in f_\theta^\alpha(A, \mathscr{M})$.\\

But,
\begin{eqnarray*}
\displaystyle\lim_{r \rightarrow \infty} \frac{1}{h_r^\alpha} \displaystyle\sum_{i \in I_r} M_i\left( \frac{\vert A_i(x)-L \vert}{\rho^{(i)}}\right)
&=&
\displaystyle\lim_{r \rightarrow \infty} \frac{1}{h_r^\alpha} M_{n_r}\left(\frac{\vert \nu_r- L\vert}{\rho^{(i)}}\right)\\
&\geq &
\displaystyle\lim_{r \rightarrow \infty} \frac{1}{h_r^\alpha} h_r^\alpha\\
&=&1
\end{eqnarray*}

So, $x \in w_\theta^\alpha(A, \mathscr{M})$.

\end{document}